\documentclass[10pt]{article}
\usepackage{amsfonts}
\usepackage{latexsym}
\usepackage{epsfig}
\usepackage{theorem}
\usepackage{amsmath}
\usepackage{amssymb}
\usepackage{graphics}
\usepackage{epsfig}
\usepackage{subfigure}
\usepackage{float}

\newenvironment{proof}{
\noindent {\bf Proof.\hskip 3mm}}%
{\mbox{}\hfill\rule{0.5em}{0.809em}\par}
\newtheorem{theorem}{Theorem}[section]
\newtheorem{lemma}[theorem]{Lemma}
\newtheorem{corollary}[theorem]{Corollary}

\newtheorem{conjecture}{\bf Conjecture }

\begin{document}

\title{\LARGE Neighbor product distinguishing total colorings of corona of subcubic graphs
\thanks{This work was supported by the National Natural Science Foundation of China (Grant No. NSFC12001332).
It was also supported by China Postdoctoral Science Foundation
Funded Project (Grant No.2014M561909); the Nature Science
Foundation of Shandong Province of China (Grant No. ZR2014AM028,
ZR2017BA009)}}
\author{Aijun Dong \quad \quad  Wenwen Zhang\\
{\small  School of Data and Computer Science,}\\ {\small Shandong
Women's University, Jinan, 250300, China}}

\date{}

\maketitle \vspace{-0.2cm}
\begin{abstract}

A proper $[k]$-total coloring $c$ of a graph $G$ is a mapping $c$
from $V(G)\bigcup E(G)$ to $[k]=\{1,2,\cdots,k\}$ such that
$c(x)\neq c(y)$ for which $x$, $y\in V(G)\bigcup E(G)$ and $x$ is
adjacent to or incident with $y$. Let $\prod(v)$ denote the
product of $c(v)$ and the colors on all the edges incident with
$v$. For each edge $uv\in E(G)$, if $\prod(u)\neq \prod(v)$, then
the coloring $c$ is called a neighbor product distinguishing total
coloring of $G$. By $\chi''_{\prod}(G)$, we denote the minimal
value of $k$ in such a coloring of $G$. In 2015, Li et al.
conjectured that $\chi''_{\prod}(G)\leq\Delta(G)+3$. In this
paper, for the corona graph of two arbitrary subcubic graphs, we
give an algorithm of a neighbor product distinguishing total
coloring and confirm the conjecture put forward by Li et al.
\end{abstract}

{\bf Keywords:} total coloring; corona graph; neighbor product
distinguishing coloring

\textbf{MSC(2010): 05C15}


\vspace{5mm}

\section{Introduction}

All graphs considered in this paper are simple and undirected. Let
$G=(V,E)$ be a graph. $V(G)$, $E(G)$, $\Delta(G)$ and $\delta(G)$
were used to denote the vertex set, edge set, maximum degree and
minimum degree of $G$, respectively. Let $N_G(u)$ be the set of
neighbors of $u$ in the graph $G$. We use $n_g$ and $n_h$ to
denote the number of vertices in graphs $G$ and $H$, respectively.
The notation and terminology used but undefined here can be found
in~\cite{bondy}.

Let $k$ be a positive integer. A total $k$-coloring $c$ is a
mapping $c: E(G)\cup V(G)\rightarrow [k]$, such that for any two
elements $x,y \in E(G)\cup V(G)$, if $xy\in E(G)$ or $x$ is
incident with $y$, then $c(x)\neq c(y)$. For each vertex $u\in
V(G)$, let $\prod(u)$ (resp. $S(u)$) denote the product (resp.
set) of colors on $u$ and the edges which are incident with $u$.
If $\prod(v)\neq \prod(u)$ (resp. $S(u)\neq S(v)$) for each edge
$uv\in E(G)$, then it is called $a$ $neighbor$ $product$
$distinguishing$ $total$ $coloring$ (resp. $neighbor$ $vertex$
$distinguishing$ $total$ $coloring$) of $G$. For convenience, they
are abbreviated as $NPDTC$ and $NVDTC$, respectively. The
chromatic number of $NPDTC$ (resp. $NVDTC$) $\chi''_{\prod}(G)$
(resp. $tndi(G)$) is the smallest integer number $k$ such that $G$
admits a neighbor product (resp. vertex) distinguishing
$[k]$-total coloring. Clearly, $tndi(G)\leq \chi''_{\prod}(G)$.

In 2005, Zhang et al.~\cite{zhangzhongfu1} determined the neighbor
set distinguishing total coloring index for graphs which are
cliques, paths, cycles, fans, wheels, stars, complete graphs etc.
and made the following conjecture.

\begin{conjecture}\label{conj1} Let $G$ be a connected graph with at least two vertices,
then $tndi(G)\leq \Delta(G)+3$.
\end{conjecture}

Wang, Chen~\cite{chenxiangen1, wanghaiying1} and Lu et
al.~\cite{luyou2} relayed to verify the Conjecture~\ref{conj1} for
graphs with $\Delta\leq4$. The Conjecture~\ref{conj1} was
confirmed for some special graphs such as $1$-tree, some sparse
graph, $2$-degenerate graph and line and splitting graph of some
graphs
in~\cite{wanghaiying2},~\cite{wangweifan2},~\cite{miaozhengke1}
and~\cite{thirusangu}, respectively. Outer planar graph,
series-parallel graphs were studied in~\cite{wangweifan1,
wangweifan3}. Chang et al.~\cite{chang1} confirmed the
Conjecture~\ref{conj1} for planar graph with $\Delta\geq8$. More
related results can be seen in ~\cite{huangdanjun, wangweifan4,
chengxiaohan1, hujie}.

Recently, Li et al.~\cite{litong} introduced the notation of NPDTC
and proposed the following conjecture.

\begin{conjecture}\label{conj2} If $G$ is a graph with at least two vertices,
then $\chi''_{\prod}(G)\leq\Delta(G)+3$.
\end{conjecture}

Li et al. prove that the conjecture holds for complete graphs,
cycles, trees, bipartite graphs, subcubic graphs and $K_4$-minor
free graphs. In 2017, Dong et al. confirmed the
Conjecture~\ref{conj2} for sparse graph $G$ with bounded maximum
degree~\cite{dong1}. Recently, the Conjecture~\ref{conj2} was
confirmed for $2$-degenerate graph in~\cite{zhuyu}.

In this paper, we discuss the NPDTC of graph products which was
introduced by Frucht and Harary in 1970~\cite{Frucht}. Given two
simple graphs $G$ and $H$, the $corona$ $product$ of $G$ and $H$
is abbreviated as $G\circ H$ which is obtained by taking one copy
of $G$, $|V(G)|$ copies of $H$ ($H_1$, $H_2$,$\ldots$,$H_{n_g}$),
and making every vertex $v_j$ of $G$ adjacent to every vertex
$u^j_i$ (the copy of $u_i$ in $H_j$) of $H_j$ to get edges
$v_ju^j_i$ ($1\leq j\leq n_g$, $1\leq i\leq n_h$). The edge
$v_ju^j_i$ is called a $corona$ $edge$. Obviously, there are
$n_g\cdot n_h$ corona edges in $G\circ H$ in total. For
convenience, let $d_G(v)$, $d_H(v)$ and $d(v)$ denote the degree
of a vertex $v$ in $G$, $H$ and $G\circ H$, respectively.
Obviously, the maximum degree of the corona graph $G\circ H$ is
$\Delta(G\circ H)=\Delta(G)+n_h$. A corona graph $K_3\circ K_4$
can be seen in Figure 1.

\begin{figure}[htbp]
  \begin{center}
    \includegraphics[width=5cm]{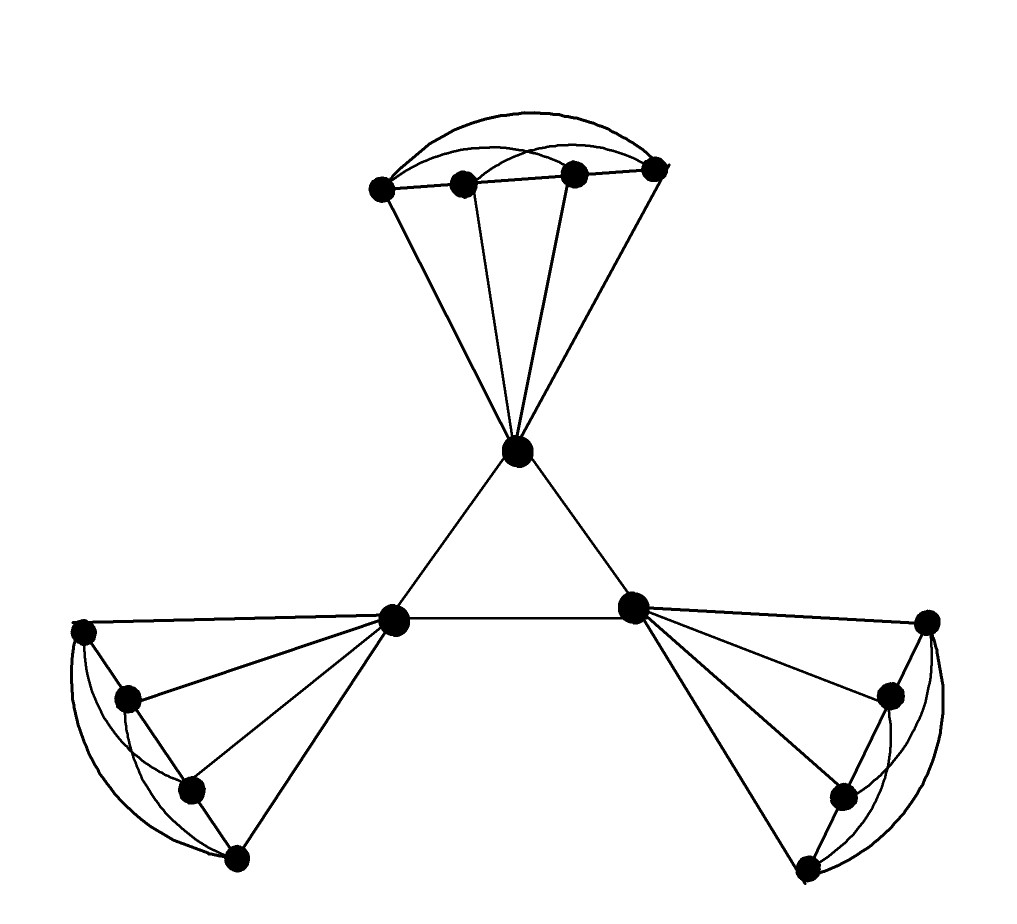}\\
   {Figure 1}
  \end{center}
\end{figure}

The concept of the corona product has some applications in
chemistry for representing chemical compounds~\cite{johnson}.
Other applications of this concept include navigation of robots in
networks~\cite{khuller}. Some theoretical results can be seen
in~\cite{daouya,mohan}. In this paper, we focus on the corona
graph of two arbitrary subcubic graphs and get the following
result.

\begin{theorem}\label{theorem1} Let $G$, $H$ be two arbitrary subcubic graphs.
Then $\chi''_{\prod}(G\circ H)\leq\Delta(G\circ H)+3$.
\end{theorem}

Since the neighbor product distinguishing $k$-total coloring is a
generation of the neighbor vertex distinguishing $k$-total
coloring, we have the following corollary.

\begin{corollary}\label{corollary1} Let $G$,$H$ be two arbitrary subcubic graphs.
Then $tndi(G\circ H)\leq\Delta(G\circ H)+3$.
\end{corollary}

\section{Preliminaries}
Let us recall some known results that will be useful for the
forthcoming proof.

\begin{lemma}[\cite{pilsniak}]\label{lm1} If $G$ is a subcubic graph,
then $\chi''_{\prod}(G)\leq \Delta(G)+3$.
\end{lemma}

\begin{lemma}[\cite{vizing}]\label{lm2} Let $G$ be a graph. Then
$\Delta(G)\leq \chi'(G)\leq \Delta(G)+1$.
\end{lemma}

\section{Proof of Theorem~\ref{theorem1}}

\begin{proof}Since $\Delta(G)\leq3$, then $G$ is a subcubic graph. We have
$\chi''_{\prod}(G)\leq\Delta(G)+3$ by Lemma~\ref{lm1}. By
Lemma~\ref{lm2}, we have $\chi'(H)\leq\Delta(H)+1$. For
convenience, we use $c''$ and $c'$ to denote a neighbor product
distinguishing $(\Delta(G)+3)$-total coloring of $G$ and a
$(\Delta(H)+1)$-edge coloring of $H$, respectively. For each
vertex $v\in V(G)$, we use $\prod_{c''}(v)$ to denote the product
of colors on the vertex $v$ and the edges which are incident with
$v$ in the coloring $c''$ of $G$, i.e.
$\prod_{c''}(v)=c''(v)\cdot\prod_{v\in e}(c''(e))$. Clearly, for
each edge $v'v''\in E(G)$, we have
$\prod_{c''}(v')\neq\prod_{c''}(v'')$. For each vertex $u\in
V(H)$, by $\prod_{c'}(u)$ (resp. $S_{c'}(u)$), we denote the
product (resp. set) of colors on the edges which are incident with
$u$ in the coloring $c'$ of $H$. Without loss of generality, we
assume that
$\prod_{c'}(u_1)\leq\prod_{c'}(u_2)\leq\ldots\leq\prod_{c'}(u_{n_h})$.
We assume $H_j$ (the $j^{th}$ copy of $H$) has the same coloring
$c'$ as $H$ for $1\leq j\leq n_g$. Obviously, we also have
$\prod_{c'}(u^j_1)\leq\prod_{c'}(u^j_2)\leq\ldots\leq\prod_{c'}(u^j_{n_h})$.
In the following, we will give colors to the corona edges
$v_ju^j_i$ and the vertices $u^j_i$ for $1\leq j\leq n_g$, $1\leq
i\leq n_h$ to get a neighbor product distinguishing
$(\Delta(G\circ H)+3)$-total coloring of $G\circ H$.

Now, we divide the proof into the following cases.

\textbf{Case $1.$ $\Delta(G)=1$}.

Let $G=v_1v_2$. Clearly, $\chi''_{\prod}(G)=3$. Let
$S=\{1,2,3,4\}$. Since $|S_{c'}(u_1)|\leq3$, then
$|(S-S_{c'}(u_1))|\geq1$.

\textbf{Case $1.1.$ $4\notin(S-S_{c'}(u_1))$}.

Since $|(S-S_{c'}(u_1))|\geq1$, there exists at least one color
$\beta\in(S-S_{c'}(u_1))$ such that $\beta\in\{1,2,3\}$. Now, let
$c(u_1^j)=c(v_1v_2)=\beta$, $c(v_1u_1^1)=c(v_2u_1^2)=5$, and use
the other two colors in $\{1,2,3\}$ to color $v_1$, $v_2$,
respectively. Let $c(u_i^j)=i+3$, $c(v_ju_i^j)=i+4$ for each
$1\leq j\leq2$ and $2\leq i\leq n_h$. So far, we get a proper
$(\Delta(G)+n_h+3)$-total coloring $c$ of $G\circ H$ such that
$c(e)=c'(e)$ for each edge $e\in E(H)$. Clearly,
$\prod_c(v_1)=c(v_1)\cdot\beta\cdot5\cdot6\cdot\ldots\cdot(n_h+3)\cdot(n_h+4)$,
$\prod_c(v_2)=c(v_2)\cdot\beta\cdot5\cdot6\cdot\ldots\cdot(n_h+3)\cdot(n_h+4)$,
and $\prod_c(u_1^j)=\prod_{c'}(u^j_1)\cdot\beta\cdot5$ for $1\leq
j\leq2$. For each vertex $u^j_i\in V(H_j)$ where $1\leq j\leq 2$
and $2\leq i\leq n_h$, we have
$\prod_c(u^j_i)=\prod_{c'}(u^j_i)\cdot(i+3)\cdot(i+4)$.

In the following, we will show the coloring $c$ is a neighbor
product distinguishing total coloring of $G\circ H$.

First, we consider the vertices of $G$. Clearly,
$\frac{\prod_c(v_1)}{\prod_c(v_2)}=\frac{c(v_1)}{c(v_2)}\neq1$. So
we have $\prod_c(v_1)\neq\prod_c(v_2)$.

Second, we will show $\prod_c(u^j_{i_1})\neq\prod_c(u^j_{i_2})$
for any two adjacent vertices $u^{j}_{i_1}$ and $u^j_{i_2}$ of
$H_j$. Since $\beta\leq3$, then
$\prod_c(u_1^j)=\prod_{c'}(u^j_{1})\cdot\beta\cdot5<\prod_{c'}(u^j_{i})\cdot(i+3)\cdot(i+4)$
for $2\leq i\leq n_h$. For each $2\leq i_1<i_2\leq n_h$, since
$\prod_{c'}(u^j_{i_1})\leq\prod_{c'}(u^j_{i_2})$, then
$\prod_{c'}(u^j_{i_1})\cdot(i_1+3)\cdot(i_1+4)<
\prod_{c'}(u^j_{i_2})\cdot(i_2+3)\cdot(i_2+4)$, i.e.
$\prod_c(u^j_{i_1})<\prod_c(u^j_{i_2})$.

At last, we consider two adjacent vertices $u^j_i$ and $v_j$ for
$1\leq j\leq 2$, $1\leq i\leq n_h$. For each $1\leq j\leq2$, we
have
$\prod_c(v_j)=\prod_{c''}(v_j)\cdot5\cdot6\cdot\ldots\cdot(n_h+3)\cdot(n_h+4)\geq
1\cdot2\cdot5\cdot6\cdot\ldots\cdot(n_h+3)\cdot(n_h+4)$.
Furthermore,
$\prod_c(u^j_1)<\prod_c(u^j_2)<\ldots<\prod_c(u^j_{n_h})$, and
$\prod_c(u^j_{n_h})=\prod_{c'}(u^j_{n_h})\cdot(n_h+3)\cdot(n_h+4)$.
Since $4\in S_{c'}(u_1)$, then $\chi'(H)=4$. Obviously, $n_h>4$.
We have
$\prod_c(u^j_{n_h})\leq2\cdot3\cdot4\cdot(n_h+3)\cdot(n_h+4)$.
Since $\prod_c(v_j)\geq
1\cdot2\cdot5\cdot6\cdot\ldots\cdot(n_h+3)\cdot(n_h+4)$, then
$\frac{\prod_c(v_j)}{\prod_c(u^j_{n_h})}\geq\frac{5}{2}$. So we
have $\Sigma_c(v_j)>\Sigma_c(v^j_i)$ for each $1\leq j\leq2$ and
$1\leq i\leq n_h$.

\textbf{Case $1.2.$ $4\in(S-S_{c'}(u_1))$}.

Let $c(u_i^j)=i+3$, $c(v_ju_i^j)=i+4$ for each $1\leq j\leq2$ and
$1\leq i\leq n_h$. So far, we get a proper
$(\Delta(G)+n_h+3)$-total coloring $c$ of $G\circ H$ such that
$c(e)=c'(e)$ for each edge $e\in E(H)$ and $c(x)=c''(x)$ for each
$x\in V(G)\bigcup E(G)$. For each vertex $v\in V(G)$, we have
$\prod_c(v)=\prod_{c''}(v)\cdot5\cdot6\cdot\ldots\cdot(n_h+3)\cdot(n_h+4)$.
For each vertex $u^j_i\in V(H_j)$ where $1\leq j\leq 2$ and $1\leq
i\leq n_h$, we have
$\prod_c(u^j_i)=\prod_{c'}(u^j_i)\cdot(i+3)\cdot(i+4)$.

In the following, we will show the coloring $c$ is a neighbor
product distinguishing total coloring of $G\circ H$.

First, we consider the vertices of $G$. Since
$\prod_{c''}(v_1)\neq\prod_{c''}(v_2)$, clearly, we have
$\prod_{c''}(v_1)\cdot5\cdot6\cdot\ldots\cdot(n_h+3)\cdot(n_h+4)\neq
\prod_{c''}(v_2)\cdot5\cdot6\cdot\ldots\cdot(n_h+3)\cdot(n_h+4)$,
i.e. $\prod_c(v_1)\neq\prod_c(v_2)$.

Second, we will show $\prod_c(u^j_{i_1})\neq\prod_c(u^j_{i_2})$
for any two adjacent vertices $u^{j}_{i_1}$ and $u^j_{i_2}$ of
$H_j$. Without loss of generality, we assume $i_1<i_2$. Since
$\prod_{c'}(u^j_{i_1})\leq\prod_{c'}(u^j_{i_2})$, then
$\prod_{c'}(u^j_{i_1})\cdot(i_1+3)\cdot(i_1+4)<
\prod_{c'}(u^j_{i_2})\cdot(i_2+3)\cdot(i_2+4)$, i.e.
$\prod_c(u^j_{i_1})<\prod_c(u^j_{i_2})$.

At last, we consider two adjacent vertices $u^j_i$ and $v_j$ for
$1\leq j\leq 2$, $1\leq i\leq n_h$. For each $1\leq j\leq2$, we
have
$\prod_c(v_j)=\prod_{c''}(v_j)\cdot5\cdot6\cdot\ldots\cdot(n_h+3)\cdot(n_h+4)\geq
1\cdot2\cdot5\cdot6\cdot\ldots\cdot(n_h+3)\cdot(n_h+4)$.
Furthermore,
$\prod_c(u^j_1)<\prod_c(u^j_2)<\ldots<\prod_c(u^j_{n_h})$, and
$\prod_c(u^j_{n_h})=\prod_{c'}(u^j_{n_h})\cdot(n_h+3)\cdot(n_h+4)$.

If $n_h=2$, then $\chi'(H)=1$. We have
$\prod_c(u^j_2)=1\cdot5\cdot6=30$. Since $\prod_c(v_j)\geq
1\cdot2\cdot5\cdot6=60$, we have $\Sigma_c(v_j)>\Sigma_c(u^j_i)$
for each $1\leq j, i\leq2$.

If $n_h=3$, then $\chi'(H)\leq3$. We have
$\prod_c(u^j_3)\leq2\cdot3\cdot6\cdot7=252$. Since
$\prod_c(v_j)\geq 1\cdot2\cdot5\cdot6\cdot7=420$, we have
$\Sigma_c(v_j)>\Sigma_c(u^j_i)$ for each $1\leq j\leq 2$, $1\leq
i\leq3$.

Otherwise, i.e. $n_h\geq4$. Clearly, $2\leq\chi'(H)\leq4$. We have
$\prod_c(u^j_{n_h})\leq2\cdot3\cdot4\cdot\ldots\cdot(n_h+3)\cdot(n_h+4)$.
Since $\prod_c(v_j)\geq
1\cdot2\cdot5\cdot6\cdot\ldots\cdot(n_h+3)\cdot(n_h+4)$, clearly,
$\frac{\prod_c(v_j)}{\prod_c(u^j_{n_h})}\geq\frac{5}{2}$. So we
have $\Sigma_c(v_j)>\Sigma_c(v^j_i)$ for each $1\leq j\leq 2$ and
$1\leq i\leq n_h$.

\textbf{Case $2.$ $2\leq\Delta(G)\leq3$}.

Clearly, $\chi''_{\prod}(G)\leq6$. Let $S=\{1,2,3,4,5\}$. Since
$|S_{c'}(u_1)\bigcup c''(v)|\leq4$, then $|S-(S_{c'}(u_1)\bigcup
c''(v))|\geq1$. There exists at least one color $\alpha\in
(S-(S_{c'}(u_1)\bigcup c''(v)))$. Let $c(u_1^j)=\alpha$,
$c(v_ju_1^j)=\Delta(G)+4$, and let $c(u_i^j)=\Delta(G)+i+2$,
$c(v_ju_i^j)=\Delta(G)+i+3$ for each $1\leq j\leq n_g$ and $2\leq
i\leq n_h$. So far, we get a proper $(\Delta(G)+n_h+3)$-total
coloring $c$ of $G\circ H$ such that $c(e)=c'(e)$ for each edge
$e\in E(H)$ and $c(x)=c''(x)$ for each $x\in V(G)\bigcup E(G)$.
For each vertex $v\in V(G)$, we have
$\prod_c(v)=\prod_{c''}(v)\cdot(\Delta(G)+4)\cdot(\Delta(G)+5)\cdot\ldots\cdot(\Delta(G)+n_h+3)$,
$\prod_c(u^j_1)=\prod_{c'}(u^j_1)\cdot\alpha\cdot(\Delta(G)+4)$.
For each vertex $u^j_i\in V(H_j)$ where $1\leq j\leq n_g$ and
$2\leq i\leq n_h$, we have
$\prod_c(u^j_i)=\prod_{c'}(u^j_i)\cdot(\Delta(G)+i+2)\cdot(\Delta(G)+i+3)$.

In the following, we will show the coloring $c$ is a neighbor
product distinguishing total coloring of $G\circ H$.

First, we consider two adjacent vertices $v'$ and $v''$ in $G$.
Since $\prod_{c''}(v')\neq\prod_{c''}(v'')$, it is obvious that
$\prod_{c''}(v')\cdot(\Delta(G)+4)\cdot(\Delta(G)+5)\cdot\ldots\cdot(\Delta(G)+n_h+3)\neq
\prod_{c''}(v'')\cdot(\Delta(G)+4)\cdot(\Delta(G)+5)\cdot\ldots\cdot(\Delta(G)+n_h+3)$,
i.e. $\prod_c(v')\neq\prod_c(v'')$.

Second, we will show $\prod_c(u^j_{i_1})\neq\prod_c(u^j_{i_2})$
for any two adjacent vertices $u^{j}_{i_1}$ and $u^j_{i_2}$ of
$H_j$. Since $\alpha\leq5$, then
$\prod_c(u_1^j)=\prod_{c'}(u^j_{1})\cdot\alpha\cdot(\Delta(G)+4)<\prod_{c'}(u^j_{i})\cdot(\Delta(G)+i+2)\cdot(\Delta(G)+i+3)$
for $2\leq i\leq n_h$. For each $2\leq i_1<i_2\leq n_h$, since
$\prod_{c'}(u^j_{i_1})\leq\prod_{c'}(u^j_{i_2})$, it is clear that
$\prod_{c'}(u^j_{i_1})\cdot(\Delta(G)+i_1+2)\cdot(\Delta(G)+i_1+3)<
\prod_{c'}(u^j_{i_2})\cdot(\Delta(G)+i_2+2)\cdot(\Delta(G)+i_2+3)$,
i.e. $\prod_c(u^j_{i_1})<\prod_c(u^j_{i_2})$.

At last, we consider two adjacent vertices $u^j_i$ and $v_j$ for
$1\leq j\leq n_g$, $1\leq i\leq n_h$. For each $1\leq j\leq n_g$,
we have
$\prod_c(v_j)=\prod_{c''}(v_j)\cdot(\Delta(G)+4)\cdot(\Delta(G)+5)\cdot\ldots\cdot
(\Delta(G)+n_h+3)\geq1\cdot2\cdot(\Delta(G)+4)\cdot(\Delta(G)+5)\cdot\ldots\cdot
(\Delta(G)+n_h+3)$. Furthermore,
$\prod_c(u^j_1)<\prod_c(u^j_2)<\ldots<\prod_c(u^j_{n_h})$, and
$\prod_c(u^j_{n_h})=\prod_{c'}(u^j_{n_h})\cdot(\Delta(G)+n_h+2)\cdot(\Delta(G)+n_h+3)$.

If $n_h=2$, then $\chi'(H)=1$. We have
$\prod_c(u^j_2)=1\cdot(\Delta(G)+4)\cdot(\Delta(G)+5)$. Since
$\prod_c(v_j)\geq 1\cdot2\cdot(\Delta(G)+4)\cdot(\Delta(G)+5)$,
then $\frac{\prod_c(v_j)}{\prod_c(u^j_2)}\geq2$. Clearly, we have
$\Sigma_c(v_j)>\Sigma_c(u^j_i)$ for each $1\leq j\leq n_g$, $1\leq
i\leq2$.

If $n_h=3$, then $\Delta(H)\leq2$, $\chi'(H)\leq3$. We have
$\prod_c(u^j_3)\leq2\cdot3\cdot(\Delta(G)+5)\cdot(\Delta(G)+6)$.
Since $\prod_c(v_j)\geq
1\cdot2\cdot(\Delta(G)+4)\cdot(\Delta(G)+5)(\Delta(G)+6)$, then
$\frac{\prod_c(v_j)}{\prod_c(u^j_3)}\geq\frac{5}{3}$. Clearly, we
have $\Sigma_c(v_j)>\Sigma_c(u^j_i)$ for each $1\leq j\leq n_g$,
$1\leq i\leq3$.

Otherwise, i.e. $n_h\geq4$. Clearly, $\chi'(H)\leq4$. We have
$\prod_c(u^j_{n_h})\leq2\cdot3\cdot4\cdot(\Delta(G)+n_h+2)\cdot(\Delta(G)+n_h+3)$.
Since $\prod_c(v_j)\geq
1\cdot2\cdot(\Delta(G)+4)\cdot(\Delta(G)+5)\cdot\ldots\cdot(\Delta(G)+n_h+2)\cdot(\Delta(G)+n_h+3)$,
clearly, $\frac{\prod_c(v_j)}{\prod_c(u^j_{n_h})}\geq\frac{5}{2}$.
So we have $\Sigma_c(v_j)>\Sigma_c(v^j_i)$ for each $1\leq j\leq
n_g$ and $1\leq i\leq n_h$.

\end{proof}



\end{document}